\documentclass[12pt]{article}
\pdfoutput=1
\usepackage{latexsym}
\usepackage{amsmath}
\usepackage{amssymb}
\usepackage{amsthm}
\usepackage{stmaryrd}
\usepackage{epsfig}
\usepackage{graphicx}
\usepackage{enumerate}
\usepackage{subfigure}
\usepackage{url}
\usepackage{tikz}
\usepackage{pgf}

\newtheorem{theorem}{Theorem}
\newtheorem{thm}{Theorem}
\newtheorem{corollary}[theorem]{Corollary}
\newtheorem{lemma}[theorem]{Lemma}

\newtheorem{remark}[theorem]{Remark}

\newtheorem{definition}[theorem]{Definition}

\newtheorem{assumption}[theorem]{Assumption}

\newcommand{\be}{\begin{enumerate}}
\newcommand{\ee}{\end{enumerate}}

\newcommand{\R}{\mathbb{R}}

\def\trace{{\bf Tr}}

\advance \topmargin by -\headheight
\advance \topmargin by -\headsep
\evensidemargin \oddsidemargin
\usepackage[pdftex,
    pdftitle={Scheduling Kalman Filters in Continuous Time},
    pdfauthor={Jerome Le Ny},
    colorlinks,linkcolor={blue},citecolor={blue},urlcolor={red},
]
{hyperref}

\title{\LARGE \bf
Scheduling Kalman Filters in Continuous Time
}

\author{ \parbox{5 in}{\centering Jerome Le Ny, Eric Feron and Munther Dahleh}
         \thanks{J. Le Ny is with the Department of Electrical and Systems Engineering, University of Pennsylvania, Philadelphia, PA 19104, USA {\tt\small jeromel@seas.upenn.edu.} M. Dahleh is with the Laboratory for Information and Decision Systems, Massachusetts Institute of Technology, Cambridge, MA 02139-4307, USA {\tt\small dahleh@mit.edu.} E. Feron is with the School of Aerospace Engineering, Georgia Tech, Atlanta, GA 30332, USA {\tt\small eric.feron@aerospace.gatech.edu.}}\\
}

\begin{document}
\maketitle

\begin{abstract}
A set of $N$ independent Gaussian  linear time invariant systems is observed by $M$ sensors whose task is to provide the best possible steady-state causal minimum mean square estimate of the state of the systems, in addition to minimizing a steady-state measurement cost. The sensors can switch between systems instantaneously, and there are additional resource constraints, for example on the number of sensors which can observe a given system simultaneously. 
We first derive a tractable relaxation of the problem, which provides a bound on the achievable performance. This bound can be computed by solving a convex program involving linear matrix inequalities. Exploiting the additional structure of the sites evolving independently, we can decompose this program into coupled smaller dimensional problems. In the scalar case with identical sensors, we give an analytical expression of an index policy proposed in a more general context by Whittle. In the general case, we develop open-loop periodic switching policies  whose performance matches the bound arbitrarily closely.
\end{abstract}


\section{Introduction}

Advances in sensor networks and the development of unmanned vehicle systems for intelligence, reconnaissance and surveillance missions require the development of data fusion schemes that can handle measurements originating from a large number of sensors observing a large number of targets, see e.g. \cite{Chakrabarty02_SN, Cortes04_coverage}. These problems have a long history \cite{Meier67_TAC}, and can be used to formulate static sensor scheduling problems as well as trajectory optimization problems for mobile sensors \cite{tiwari06_thesis, williams07_thesis}.

In this paper, we consider $M$ mobile sensors tracking the state of $N$ sites or targets in continuous time. We assume that the sites can be described by $N$ plants with independent linear time invariant dynamics,
 \[
\dot x_i=A_i x_i + B_i u_i + w_i, \quad x_i(0)=x_{i,0}, \quad i=1,\ldots,N.
\]
We assume that the plant controls $u_i(t)$ are deterministic and \emph{known} for $t \geq 0$. Each driving noise $w_i(t)$ is a stationary white Gaussian noise process with zero mean and known power spectral density matrix $W_i$: 
\[
\text{Cov}(w_i(t),w_i(t'))=W_i \, \delta(t-t').
\]
The initial conditions are random variables with known mean $\bar x_{i,0}$ and covariance matrices $\Sigma_{i,0}$. By independent 
systems we mean that the noise processes of the different plants are independent, as well as the initial conditions $x_{i,0}$. Moreover the initial conditions are assumed independent of the noise processes. We shall assume in addition that 
\begin{assumption}	\label{assumption: pos-def initial condition}
The matrices $\Sigma_{i,0}$ are positive definite for all $i \in \{1,\ldots,N\}$.
\end{assumption}
This can be achieved by adding an arbitrarily small multiple of the identity matrix to a potentially non invertible matrix $\Sigma_{i,0}$. This assumption is needed in our discussion to be able to use the information filter later on and to use a technical theorem on the convergence of the solutions of a periodic Riccati equation in section \ref{section: switching policies performance - KF}.

We assume that we have at our disposal $M$ sensors to observe the N plants. If sensor $j$ is used to observe plant $i$, we obtain measurements
\[
y_{ij}=C_{ij} x_i+v_{ij}.
\]
Here $v_{ij}$ is a stationary white Gaussian noise process with power spectral density matrix $V_{ij}$, assumed positive definite. Also, $v_{ij}$ is independent of the other measurement noises, process noises, and initial states. Finally, to guarantee convergence of the filters later on, we assume throughout that

\begin{assumption}	\label{assumption: observability} 
For all $i \in \{1,\ldots,N\}$, there exists a set of indices $j_1, j_2,\ldots,j_{n_i} \in  \{1,\ldots,M\}$ such that the pair $(A_i,\tilde C_{i})$ is detectable, where $\tilde C_{i}=[C^T_{ij_1}, \ldots, C^T_{i j_{n_i}}]^T$.
\end{assumption}
\begin{assumption}	\label{assumption: controllability}  
For all $i \in \{1,\ldots,N\}$, the pair $(A_i,W_i^{1/2})$ is controllable.
\end{assumption}

Let us define
\[
\pi_{ij}(t)=\begin{cases}
1 \text{ if plant } i \text{ is observed at time $t$ by sensor } j\\
0 \text{ otherwise.}
\end{cases}
\]

We assume that each sensor can observe at most one system at each instant, hence we have the constraint
\begin{equation}	\label{eq: pathwise resource constraint}
\sum_{i=1}^N \pi_{ij}(t) \leq 1, \quad \forall t, \; \; j=1,\ldots,M.	
\end{equation}
If instead sensor $j$ is required to be always operated, constraint (\ref{eq: pathwise resource constraint}) should simply be changed to 
\begin{equation}	\label{eq: pathwise resource constraint bis}
\sum_{i=1}^N \pi_{ij}(t) = 1.
\end{equation}
The equality constraint is useful in scenarios involving sensors mounted on unmanned vehicles for example, where it might not be possible to withdraw a vehicle from operation during the mission. The performance will be worse in general than with an inequality constraint once we introduce operation costs. 

We also add the following constraint, similar to the one used by Athans \cite{Athans72_sensor}. We suppose that each system can be observed by at most one sensor at each instant, so we have
\begin{equation}	\label{eq: pathwise resource constraint II}
\sum_{j=1}^M \pi_{ij}(t) \leq 1, \quad \forall t, \; \; i=1,\ldots,N.
\end{equation}
Similarly if system $i$ must always be observed by some sensor, constraint (\ref{eq: pathwise resource constraint II}) can be changed to an equality constraint 
\begin{equation}	\label{eq: pathwise resource constraint II bis}
\sum_{j=1}^M \pi_{ij}(t) = 1.
\end{equation}
Note that a sensor in our discussion can correspond to a combination of several physical sensors, and so the constraints above can capture seemingly more general problems where we allow for example more that one simultaneous measurements per system. Using (\ref{eq: pathwise resource constraint II bis}) we could also impose a constraint on the total number of allowed observations at each time. Indeed, consider a constraint of the form
\[
\sum_{i=1}^N \sum_{j=1}^M \pi_{ij} (t) \leq p, \quad \text{for some positive integer } p. 
\]
This constraint means that $M-p$ sensors are required to be idle at each time. So we can create $M-p$ ``dummy" systems (we should choose simple scalar stable systems to minimize computations), and associate the constraint (\ref{eq: pathwise resource constraint II bis}) to each of them. Then we simply do not include the covariance matrix of these systems in the objective function (\ref{eq: objective, average init}) below.

We consider an infinite-horizon average cost problem. The parameters of the model are assumed known. We wish to design an observation policy $\pi(t)=\{ \pi_{ij}(t) \}$ satisfying the constraints (\ref{eq: pathwise resource constraint}),  (\ref{eq: pathwise resource constraint II}), or their equality versions, and an estimator $\hat x_\pi$ of $x$, \emph{depending at each instant only on the past and current observations produced by the observation policy}, such that the average error covariance is minimized, in addition to some observation costs. The policy $\pi$ itself can also only depend on the past observations. More precisely, we wish to minimize, subject to the constraints (\ref{eq: pathwise resource constraint}),  (\ref{eq: pathwise resource constraint II}),
\begin{equation}	\label{eq: objective, average init}
J_{avg}=\min_{\pi,\hat x_\pi} \limsup_{T \to \infty} \frac{1}{T} E \left[ \int_0^T  \sum_{i=1}^N \left ( (x_i-\hat x_{\pi,i})' T_i (x_i-\hat x_{\pi,i}) + \sum_{j=1}^M \kappa_{ij} \, \pi_{ij}(t) \right ) dt \right],
\end{equation}
where the constants $\kappa_{ij}$ are a cost paid per unit of time when plant $i$ is observed by sensor $j$. The $T_i$'s are positive semidefinite weighting matrices.

\emph{Literature Review and Contributions of this paper}. The sensor scheduling problem presented above, except for minor variations, is an infinite horizon version of the problem considered by Athans in \cite{Athans72_sensor}. See also  Meier et al. \cite{Meier67_TAC} for the discrete-time version. Athans considered the observation of only one plant. We include here several plants to show how their independent evolution property can be leveraged in the computations, using the dual decomposition method from optimization.

Discrete-time versions of this sensor selection problem have received a significant amount of attention, see e.g. \cite{ feron90_KF, oshman94_sensor, Tiwari05_dynamicCoverage, tiwari06_thesis,Gupta06_sensorSelection, LaScala06_RBTracking, shi07_estimation}. All algorithms proposed so far, except for the optimal greedy policy of \cite{LaScala06_RBTracking} in the completely symmetric case, either run in exponential time or consist of heuristics with no performance guarantee. 
We do not consider the discrete-time problem in this paper. Finite-horizon continuous-time versions of the problem, besides the presentation of Athans  \cite{Athans72_sensor}, have also been the subject of several papers \cite{skafidas98_oms, savkin00_optimalScheduling, baras89_optimalMeasurement, lee01_sensorScheduling}. The solutions proposed, usually based on optimal control techniques, also involve computational procedures that scale poorly with the dimension of the problem.

Somewhat surprisingly however, and with the exception of \cite{mourikis06_sensorScheduling}, it seems that the infinite-horizon continuous time version of the Kalman filter scheduling problem has not been considered previously. Mourikis and Roumeliotis \cite{mourikis06_sensorScheduling} consider initially also a discrete time version of the problem for a particular robotic application. However, their discrete model originates from the sampling at high rate of a continuous time system. To cope with the difficulty of determining a sensor schedule, they assume instead a model where each sensor can independently process each of the available measurements at a constant frequency, and seek the optimal measurement frequencies. In fact, they obtain these frequencies by introducing heuristically a continuous time Riccati equation, and show that the frequencies can then be computed by solving a semidefinite program. In contrast, we consider the more standard schedule-based version of the problem in continuous time, which is a priori more constraining. We show that essentially the same convex program provides in fact a \emph{lower bound} on the cost achievable by \emph{any} measurement policy. In addition, we provide additional insight into the decomposition of the computations of this program, which can be useful in the framework of \cite{mourikis06_sensorScheduling} as well.

The rest of the chapter is organized as follows. Section \ref{section: optimal estimator} briefly recalls that for a fixed policy $\pi(t)$, the optimal estimator is obtained by a type of Kalman-Bucy filter. The properties of the Kalman filter (independence of the error covariance matrix with respect to measurement values) imply that the remaining problem of finding the optimal scheduling policy $\pi$ is a deterministic control problem. 
In section \ref{section: 1D-dynamics} we treat a simplified scalar version of the problem with identical sensors as a special case of the classical ``Restless Bandit Problem'' (RBP) \cite{whittle-restless}, and provide analytical expressions for an index policy and for the elements necessary to compute efficiently a lower bound on performance, both of which were proposed in the general setting of the RBP by Whittle. Then, for the multidimensional case treated in full generality in section \ref{section: MIMO case}, we show that the lower bound on performance can be computed as a convex program involving linear matrix inequalities. This lower bound can be approached arbitrarily closely by a family of new periodically switching policies described in section \ref{section: BvN policies}. Approaching the bound with these policies is limited only by the frequency with which the sensors can actually switch between the systems. In general, our solution has much more attractive computational properties than the solutions proposed so far for the finite-horizon problem.

\section{Optimal Estimator}	\label{section: optimal estimator}

For a given observation policy $\pi(t)=\{ \pi_{ij}(t) \}_{i,j}$, the minimum variance filter is given by the Kalman-Bucy filter \cite{kalmanBucy61_filter}, see \cite{Athans72_sensor}.
The state estimates $\hat x_\pi$, where the subscript indicates the dependency on the policy $\pi$, are all updated in parallel following the stochastic differential equation
\begin{align*}
& \frac{d}{dt}{\hat x}_{\pi,i}(t) = A_i \hat x_{\pi,i}(t) + B_i(t) u_i(t) + \Sigma_{\pi,i}(t) \left( \sum_{j=1}^M \pi_{ij}(t) C^T_{ij} V_{ij}^{-1} \left(C_{ij} \hat x_{\pi,i}(t) - y_{ij}(t) \right) \right), \\
& \hat x_{\pi,i}(0)=\bar x_{i,0}.
\end{align*}
The resulting estimator is unbiased and the error covariance matrix $\Sigma_{\pi,i}(t)$ for site $i$ verifies the matrix Riccati differential equation
\begin{align}	\label{eq: Riccati ODE}
& \frac{d}{dt} \Sigma_{\pi,i}(t)= A_i \Sigma_{\pi,i}(t) + \Sigma_{\pi,i}(t) A_i^T + W_i - \Sigma_{\pi,i}(t) \left( \sum_{j=1}^M  \pi_{ij}(t) C_{ij}^T V_{ij}^{-1} C_{ij}\right) \Sigma_{\pi,i}(t), \\
& \Sigma_{\pi,i}(0)=\Sigma_{i,0}. \nonumber
\end{align}

With this result, we can reformulate the optimization of the observation policy as a \emph{deterministic} optimal control problem. Rewriting
\[
E((x_i-\hat x_i)' T_i (x_i-\hat x_i)) = \text{Tr} \, (T_i \, \Sigma_i),
\]
the problem is to compute
\begin{equation}	\label{eq: objective, average}
\min_{\pi} \limsup_{T \to \infty} \frac{1}{T} \left[ \int_0^T \sum_{i=1}^N \left( \text{Tr} \left( T_i \, \Sigma_{\pi,i}(t) \right) + \sum_{j=1}^M \kappa_{ij} \, \pi_{ij}(t) \right ) dt \right],
\end{equation}
subject to the constraints (\ref{eq: pathwise resource constraint}), (\ref{eq: pathwise resource constraint II}), or their equality versions, and the dynamics (\ref{eq: Riccati ODE}).

\section{Sites with One-Dimensional Dynamics and Identical Sensors}	\label{section: 1D-dynamics}

We assume in this section that 
\begin{enumerate}
\item the sites or targets have one-dimensional dynamics, i.e., $x_i \in \mathbb{R}, \, i=1,\ldots,N$; and,
\item all the sensors are identical, i.e., \,  $C_{ij}=C_i, V_{ij}=V_i, \kappa_{ij}=\kappa_i, \; j=1,\ldots,M$. \label{identical sensors condition}
\end{enumerate}
Because of condition \ref{identical sensors condition}, we can simplify the problem formulation introduced above so that it corresponds exactly to a special case of the Restless Bandit Problem \cite{whittle-restless}. We define
\[
\pi_{i}(t)=\begin{cases}
1 \text{ if plant } i \text{ is observed at time $t$ by a sensor}\\
0 \text{ otherwise.}
\end{cases}
\]
Since we assumed that a system can be observed by at most one sensor, the scheduling problem is interesting only in the case $M<N$. Note that a constraint (\ref{eq: pathwise resource constraint II bis}) for some system $i$ can be eliminated, by removing one available sensor, which is always measuring the system $i$. Constraints (\ref{eq: pathwise resource constraint bis}) and (\ref{eq: pathwise resource constraint II}) can then be replaced by the single constraint
\[
\sum_{i=1}^N \pi_i(t) = M, \; \forall t.
\]
This constraint means that at each period, exactly $M$ of the $N$ sites are observed. We treat this case in this section, but again the equality sign can be changed to an inequality with very little change in our discussion. 

To obtain a lower bound on the achievable performance, we relax the constraint to enforce it only on average
\begin{equation}	\label{eq: relaxed constraint scalar}
\limsup_{T \to \infty}  \frac{1}{T} \int_0^T  \sum_{i=1}^N \pi_i(t) dt = M.
\end{equation}
Then we adjoin this constraint using a (scalar) Lagrange multiplier $\lambda$ to form the Lagrangian
\[
L(\pi, \lambda) = \limsup_{T \to \infty} \frac{1}{T}  \int_0^T \sum_{i=1}^N \left[ \text{Tr} \left( T_i \, \Sigma_{\pi,i}(t) \right) + (\kappa_i+\lambda) \, \pi_{i}(t) \right] dt - \lambda M.
\]
Here $\kappa_i$ is the cost per time unit for observing site $i$. The dynamics of $\Sigma_{\pi,i}$ are now given by
\begin{align}	\label{eq: Riccati ODE - simplified}
& \frac{d}{dt} \Sigma_{\pi,i}(t)= A_i \Sigma_{\pi,i}(t) + \Sigma_{\pi,i}(t) A_i' + W_i - \pi_i(t) \, \Sigma_{\pi,i}(t) C_{i}^T V_{i}^{-1} C_{i} \Sigma_{\pi,i}(t), \\
& \Sigma_{\pi,i}(0)=\Sigma_{i,0}. \nonumber
\end{align}
Then the original optimization problem (\ref{eq: objective, average}) with the relaxed constraint (\ref{eq: relaxed constraint scalar}) can be expressed as
\[
\gamma=\inf_\pi \sup_\lambda L(\pi,\lambda) = \sup_\lambda \inf_\pi L(\pi,\lambda),
\]
where the exchange of the supremum and the infimum can be justified using a minimax theorem for constrained dynamic programming \cite{altman-constrainedMDP}. 
We are then led to consider the computation of the dual function
\[
\gamma(\lambda) = \min_\pi \limsup_{T \to \infty} \frac{1}{T} \int_0^T \sum_{i=1}^N \left[ \text{Tr} \left( T_i \Sigma_{\pi,i}(t) \right) + (\kappa_i+\lambda) \, \pi_{i}(t) \right] dt - \lambda M, 
\]
which has the important property of \emph{being separable by site}, i.e., $\gamma(\lambda)+\lambda M=\sum_{i=1}^N \gamma^i(\lambda)$, where for each site $i$ we have
\begin{equation}	\label{eq: individual site objective}
\gamma^i(\lambda)=\min_{\pi_i}  \limsup_{T \to \infty} \frac{1}{T} \int_0^T  \text{Tr} \left( T_i \, \Sigma_{\pi_i,i}(t)  \right) + (\kappa_{i} + \lambda) \, \pi_{i}(t) dt.
\end{equation}

When the dynamics of the sites are one dimensional, i.e., $\Sigma_i \in \mathbb{R}$, we can solve this optimal control problem for each site analytically, that is, we obtain an analytical expression of the dual function, which provides a lower bound on the cost for each $\lambda$. The computations are presented in paragraph (\ref{section: analytical solution}). First, we explain how these computations will also provide the elements necessary to design a scheduling policy.

\subsection{Restless Bandits}

The Restless Bandit Problem (RBP) was introduced by Whittle in \cite{whittle-restless} as a generalization of the classical Multi-Armed Bandit Problem (MABP), which was first solved by Gittins \cite{gittinsJones}. In the RBP, we have $N$ projects evolving independently, $M$ of which can be activated at each time. Projects that are active can evolve according to different dynamics than project that remain passive. In our problem, the projects correspond to the systems and their activation corresponds to taking a measurement. We describe in our particular context the index policy proposed by Whittle for the RBP, which, although suboptimal in general, generalizes the optimal policy of Gittins' in the case of the MABP.

Consider the objective (\ref{eq: individual site objective}) for system $i$. Clearly, the Lagrange multiplier $\lambda$ can be interpreted as a tax penalizing measurements of the system. As $\lambda$ increases, the passive action (i.e., not measuring) should become more attractive. For a given value of $\lambda$, let us denote $\mathcal{P}^i(\lambda)$ the set of covariance matrices $\Sigma^i$ for which the passive action is optimal. Let $\mathbb{S}^n_+$ be the set of symmetric positive semidefinite matrices. Then we say that 

\begin{definition}	\label{def: Whittle indexability}
System $i$ is indexable if and only if $\mathcal{P}^i(\lambda)$ is monotonically increasing (in the sense of set inclusion) from $\emptyset$ to $\mathbb{S}^n_+$ as $\lambda$ increases from $-\infty$ to $+\infty$. If system $i$ is indexable, we define its \emph{Whittle index} by
$
\lambda^i(\Sigma_i)=\inf \{\lambda \in \mathbb{R}: \Sigma^i \in \mathcal{P}^i(\lambda) \}.
$
\end{definition}

However natural the indexability requirement might appear, Whittle provided an example of an RBP where it is not verified. We will see in the next paragraph however that for our particular problem, at least in the scalar case, indexability of the systems is guaranteed. The idea behind the definition of the Whittle index consists in defining an intrisic ``value" for the measurement of system $i$, taking into account both the immediate and future gains. If the covariance of system $i$ is $\Sigma^i$, the Whittle index defines this value as the measurement tax (potentially negative) that should be required to make the controller indifferent between measuring and not measuring the system. Finally, if all the systems are indexable, \emph{the Whittle policy chooses at each instant to measure the $M$ systems with highest Whittle index}. There is significant experimental data and some theoretical evidence indicating that when the Whittle policy is well-defined for an RBP, its performance is often very close to optimal, see e.g. \cite{glazebrook06-RBindices, weber90_asymptoticOptiRB, LeNy_ACC08_RB}.

\subsection{Solution of the Scalar Optimal Control Problem}	\label{section: analytical solution}

We can now consider problem (\ref{eq: individual site objective}) for a single site, dropping the index $i$. We return to the scalar case $\Sigma \in \mathbb{R}$. The dynamical evolution of the variance obeys the equation
\[
\dot \Sigma = 2 A \Sigma + W - \pi \frac{C^2}{V} \Sigma^2, \quad \text{with } \pi(t) \in \{0,1\}.
\]
The Hamilton-Jacobi-Bellman (HJB) equation is
\begin{equation}	\label{eq: HJB scalar case}
\gamma(\lambda)=\min \left \{T \, \Sigma+ (2 A \Sigma + W) h'(\Sigma;\lambda), \, T \, \Sigma+\kappa+\lambda+(2 A \Sigma + W - \frac{C^2}{V} \Sigma^2) h'(\Sigma;\lambda) \right \},
\end{equation}
where $h$ is the relative value function. We will use the following notation. Consider the algebraic Riccati equation (ARE)
\[
2A x + W - \frac{C^2}{V} x^2=0.
\]
First, if $T=0$, it is clearly optimal to always observe if $\lambda+\kappa<0$ and always observe otherwise. Hence the Whittle index is $\lambda(\Sigma)=-\kappa$ for all $\Sigma \in \R_+$, and $\gamma(\lambda)=\min\{\kappa+\lambda,0\}$. So we can now assume $T>0$. If $C=0$, the solution to (\ref{eq: HJB scalar case}) is to always observe if $(\kappa+\lambda)<0$ and never observe otherwise. Hence the Whittle index is again $\lambda(\Sigma)=-\kappa$ for all $\Sigma \in \R_+$ and we get, by letting $\Sigma=-\frac{W}{2A}$ in the HJB equation for a stable system:
\[
\text{for C=0: } \gamma(\lambda)=\begin{cases}
-T \, \frac{W}{2A}+\kappa+\lambda \text{ if the system is stable ($A<0$) and } (\kappa+\lambda)<0, \\
-T \, \frac{W}{2A} \text{ if the system is stable and } (\kappa+\lambda) \geq 0, \\
+\infty \text{ otherwise ($A \geq 0$)}.
\end{cases}
\]
The third case is clear from the fact that the system is unstable and cannot be measured. So we can now assume that $C \neq 0$. Then the ARE has two roots
\[
x_1=\frac{A-\sqrt{A^2+C^2 W/V}}{C^2/V}, \quad  x_2=\frac{A+\sqrt{A^2+C^2 W/V}}{C^2/V}.
\]
By assumption \ref{assumption: controllability}, $W \neq 0$ and so $x_1$ is strictly negative and $x_2$ is strictly positive.

We can treat the case $\kappa+\lambda < 0$ immediately. Then it is obviously optimal to always observe, and we get, letting $\Sigma=x_2$ in the HJB equation:
\[
\gamma(\lambda)= T \, x_2+\kappa+\lambda.
\]
So from now on we can assume $\lambda \geq -\kappa$. 
Let us temporarily assume the following result on the form of the optimal policy. The validity of this assumption can be verified a posteriori from the formulas obtained below, using the fact that the dynamic programming equation provides a sufficient condition for optimality of a solution.

\emph{Form of the optimal policy}.
The optimal policy is a threshold policy, i.e., it observes the system for $\Sigma \geq \Sigma_{th}$ and does not observe for $\Sigma<\Sigma_{th}$, for some $\Sigma_{th} \in \mathbb{R}_+$.

We would like to obtain the value of the average cost $\gamma(\lambda)$ and of the threshold $\Sigma_{th}(\lambda)$.
Note that we already know $\Sigma_{th}(\lambda)=0$ for $\lambda \leq -\kappa$, and we have $\mathcal{P}(\lambda)=[0,\Sigma_{th}(\lambda)]$ for the passive region $\mathcal P(\lambda)$ of definition \ref{def: Whittle indexability}. Then the system is indexable if and only if $\Sigma_{th}(\lambda)$ is an increasing function of $\lambda$, and then inverting the relation $\lambda \mapsto \Sigma_{th}(\lambda)$ gives the Whittle index $\Sigma \mapsto \lambda(\Sigma)$.

\subsubsection{Case $\Sigma_{th} \leq x_2$}

In this case, we obtain as before
\begin{equation}	\label{eq: avge cost case 1}
\gamma(\lambda)= T \, x_2+\kappa+\lambda.
\end{equation}
This is intuitively clear for $\Sigma_{th}<x_2$: even when observing the system all the time, the variance still converges in finite time to a neighborhood of $x_2$. Since this neighborhood is in the active region by hypothesis, after potentially a transition period (if the variance started at a value smaller than $\Sigma_{th}$), we should always observe, and so the infinite-horizon average cost is the same as for the policy that always observes.

By continuity at the interface between the active and passive regions, we have
\begin{align*}
T \, \Sigma_{th}+(2A \Sigma_{th}+W)h'(\Sigma_{th})&=T \, \Sigma_{th}+(\kappa+\lambda)+(2A\Sigma_{th}+W-\frac{C^2}{V}\Sigma_{th}^2) h'(\Sigma_{th}) \\
\text{i.e., } \kappa+\lambda&=\frac{C^2}{V}\Sigma_{th}^2 h'(\Sigma_{th}).
\end{align*}
We have then
\begin{align}
&\frac{C^2}{V}\Sigma_{th}^2(T \, x_2+\kappa+\lambda)=\frac{C^2}{V}T\, \Sigma_{th}^3+(2A\Sigma_{th}+W)(\kappa+\lambda) \nonumber \\
&\left(-\frac{C^2}{V}\Sigma_{th}^2+2A\Sigma_{th}+W \right)(\kappa+\lambda)=\frac{C^2}{V}\Sigma_{th}^2 T\, (x_2-\Sigma_{th}) \nonumber \\
&-(\Sigma_{th}-x_2)(\Sigma_{th}-x_1)(\kappa+\lambda)=\Sigma_{th}^2 T\, (x_2-\Sigma_{th})		\quad \text{(since $C \neq 0$)} \nonumber \\
&T \, \Sigma_{th}^2-\Sigma_{th}(\kappa+\lambda)+x_1(\kappa+\lambda)=0 \nonumber \\
\text{so } \quad &\lambda(\Sigma_{th})=-\kappa+\frac{T \Sigma_{th}^2}{\Sigma_{th}-x_1}, 
\quad \Sigma_{th}(\lambda)=\frac{ \frac{\kappa+\lambda}{T}+\sqrt{( \frac{\kappa+\lambda}{T})( \frac{\kappa+\lambda}{T}-4 x_1)} }{2}. \label{eq: threshold case 1}
\end{align}
Expressions (\ref{eq: avge cost case 1}) and (\ref{eq: threshold case 1}) are valid under the condition $\Sigma_{th}(\lambda) \leq x_2.$
Note from (\ref{eq: threshold case 1}) that $\Sigma_{th} \mapsto \lambda(\Sigma_{th})$ is an increasing function and the functions $\lambda(\cdot)$ and $\Sigma_{th}(\cdot)$ are inverse of each other.

\subsubsection{Case $\Sigma_{th} > x_2$}

It turns out that in this case we must distinguish between stable and unstable systems. For a stable system ($A<0$), the Lyapunov equation 
\[
2Ax+W=0
\]
has a strictly positive solution $x_e=-\frac{W}{2A}$, with $x_e > x_2$ since $C \neq 0$.

\vspace{0.5cm}

\textbf{Stable System ($A<0$) with $\Sigma_{th} \geq x_e$}:

In this case we know that $x_e$ is in the passive region. Hence, with $\Sigma=x_e$ in the HJB equation, we get
\begin{equation}
\gamma(\lambda)= T x_e.	\label{eq: avge cost case 2}
\end{equation}
Then we have again
$
\kappa+\lambda=\frac{C^2}{V} \Sigma_{th}^2 h'(\Sigma_{th}),
$
and now 
$
T \Sigma_{th}+(2A\Sigma_{th}+W)h'(\Sigma_{th})=T x_e.
$
Hence
\begin{align}
&\frac{C^2}{V}\Sigma_{th}^2 T (x_e-\Sigma_{th})=(2A\Sigma_{th}+W)(\kappa+\lambda)=2A(\Sigma_{th}-x_e)(\kappa+\lambda) \nonumber \\
&\lambda(\Sigma_{th})=-\kappa+\frac{C^2 T \Sigma_{th}^2}{2|A|V}, \quad \Sigma_{th}(\lambda)=\frac{\sqrt{2 |A| V \frac{\lambda+\kappa}{T}}}{|C|}. 	\label{eq: threshold case 2}
\end{align}

\textbf{Stable System ($A<0$) with $x_2< \Sigma_{th} < x_e$, or non-stable system ($A \geq 0$)}:

If the system is marginally stable or unstable, we cannot define $x_e$. We can think of this case as $x_e \to \infty$ as $A \to 0_{-}$, and treat it simultaneously with the case where the system is stable and $x_2 < \Sigma_{th}< x_e$. Then $x_2$ is in the passive region, and $x_e$ is in the active region, so the prefactors of $h'(x)$ in the HJB equation do not vanish. There is no immediate relation providing the value of $\gamma(\lambda)$. We can use the smooth-fit principle to handle this case and obtain the expression of the Whittle indices, following \cite{whittle-restless}. Again the formal justification comes from using the final expressions of the value function thus obtained to verify that it indeed satisfies the HJB equation.


\begin{thm}[\cite{whittle-restless},\cite{dusonchet03_thesis}] \label{thm: 1d deterministic RBP}
Consider a continuous-time one-dimensional restless bandit project $x(t) \in \mathbb{R}$ satisfying
\[
\dot x(t)=a_k(x), \quad k=0,1,
\]
with passive and active cost rates $r_k(x), \; k=0,1$. Assume that $a_0(x)$ does not vanish in the optimal passive region, and $a_1(x)$ does not vanish in the optimal active region. Then the Whittle index is given by
\[
\lambda(x)=r_0(x)-r_1(x)+\frac{[a_1(x)-a_0(x)][a_0(x)r_1'(x)-a_1(x) r_0'(x)]}{a_0(x)a_1'(x)-a_1(x)a_0'(x)}.
\]
\end{thm}

\begin{remark}
The assumption that $a_0$ and $a_1$ do not vanish in the optimal passive and active regions respectively excludes the cases previously studied. It is missing from \cite{whittle-restless}, \cite{dusonchet03_thesis}, which therefore provide only an incomplete description of the Whittle indices for one-dimensional continuous-time deterministic projects.
\end{remark}

\begin{proof}
The derivation of the expression of the Whittle index can be found in \cite{whittle-restless}, \cite[p.53]{dusonchet03_thesis}, and is valid only under the additional assumption mentioned above.
\end{proof}

\vspace{0.1cm}

\begin{corollary}
The Whittle index for the case $x_2<\Sigma_{th}<x_e$ is given by:
\begin{equation}	\label{eq: W indices KF average}
\lambda(\Sigma_{th})=-\kappa + \frac{C^2}{2V} \frac{T \Sigma_{th}^3}{A \Sigma_{th} + W}. 
\end{equation} 
 \end{corollary}

\vspace{0.1cm}

\begin{proof}
For $x_2<\Sigma_{th}<x_e$, the assumptions of theorem \ref{thm: 1d deterministic RBP} are verified with
\begin{align*}
a_0(\Sigma)=2 A \Sigma + W, \quad & a_1(\Sigma)=2 A \Sigma + W-\frac{C^2}{V}\Sigma^2 \\
r_0(\Sigma)=T \Sigma, \quad & r_1(\Sigma)=T \Sigma+\kappa.
\end{align*}
The result follows by a straightforward calculation. Note the expression for $\lambda(\Sigma)$ indeed makes sense since we can verify that it defines an increasing function of $\Sigma$.

\end{proof}

With the value of the Whittle index, we can finish the computation of the lower bound $\gamma(\lambda)$ for the case $x_2<\Sigma_{th}<x_e$. 
Inverting the relation (\ref{eq: W indices KF average}), we obtain, for a given value of $\lambda$, the boundary $\Sigma_{th}(\lambda)$ between the passive and active regions. $\Sigma_{th}(\lambda)$ verifies the depressed cubic equation
\begin{equation}	\label{eq: cubic}
X^3 - \frac{2V(\lambda + \kappa)}{TC^2} A \, X -  \frac{2V(\lambda + \kappa)}{TC^2} W = 0.
\end{equation}
For $\lambda+\kappa \geq 0$, by Descartes' rule of signs, this polynomial has exactly one positive root, which is $\Sigma_{th}(\lambda)$. 

The HJB equation then reduces to 
\begin{align}
\gamma(\lambda)&=T \Sigma+h'(\Sigma)(2 A \Sigma+W), \quad \text{for } \Sigma < \Sigma_{th}(\lambda)  \label{eq: VF passive - average cost} \\
\gamma(\lambda)&=T \Sigma+\kappa+\lambda+h'(\Sigma)(2 A \Sigma+W- \frac{C^2}{V} \Sigma^2), \quad \text{for } \Sigma \geq \Sigma_{th}(\lambda). \label{eq: VF active - average cost}
\end{align}
Now for $x_2 < \Sigma_{th}(\lambda)<x_e$, letting $x=\Sigma_{th}(\lambda)>0$ in the HJB equation, assuming continuity of $h'$ at the boundary of the passive and active regions and eliminating $h'(\Sigma_{th}(\lambda))$, we get
\begin{align*}
&\gamma(\lambda)=T \Sigma_{th}(\lambda)+\kappa+\lambda+(\gamma-T\Sigma_{th}(\lambda)) \left(1-\frac{C^2}{V} \frac{ (\Sigma_{th}(\lambda))^2 }{2 A \Sigma_{th}(\lambda)+W} \right) \\
&(\gamma (\lambda)-T \Sigma_{th}(\lambda)) \left( \frac{C^2}{V} \frac{ (\Sigma_{th}(\lambda))^2 }{2 A \Sigma_{th}(\lambda)+W} \right) = \kappa+\lambda \\
& \gamma(\lambda)=T \Sigma_{th}(\lambda)+\frac{V (\kappa+\lambda) (2 A \Sigma_{th}(\lambda)+W)}{C^2 (\Sigma_{th}(\lambda))^2}, \quad \text{for } x_2 < \Sigma_{th}(\lambda) < x_e.
\end{align*}

\subsubsection{Summary}	\label{section: summary Whittle}

We collect the previous computations in the following theorem.

\begin{thm}
In the one-dimensional Kalman filter scheduling problem with identical sensors, the systems are indexable. For system $i$, the Whittle index $\lambda_i(\Sigma_i)$ is given as follows:
\begin{itemize}
\item Case $C=0$ or $T=0$: $\lambda_i(\Sigma_i)=-\kappa_i$, for all $\Sigma_i \in \R_+$.
\item Case $C \neq 0$ and $T \neq 0$:
\begin{equation*}
\lambda_i(\Sigma_i)=\begin{cases}
-\kappa_i+\frac{T_i \Sigma_i^2}{\Sigma_i-x_{1,i}} \; \text{ if } \Sigma_i \leq x_{2,i}, \\
-\kappa_i + \frac{C_i^2}{2V_i} \frac{T_i \Sigma_i^3}{A_i \Sigma_i + W_i} \; \text{ if } x_{2,i} < \Sigma_i < x_{e,i}, \\
-\kappa_i +\frac{T_i C_i^2 \Sigma_i^2}{2 |A_i| V_i}  \; \text{ if } x_{e,i} \leq \Sigma_i,
\end{cases}
\end{equation*}
\end{itemize}
with the convention $x_{e,i}=+\infty$ if $A_i \geq 0$. The lower bound on the achievable performance is obtained by maximizing the concave function 
\begin{equation}	\label{dual function - summary}
\gamma(\lambda)=\sum_{i=1}^N \gamma^i(\lambda) - \lambda M
\end{equation}
over $\lambda$, where the term $\gamma^i(\lambda)$ is given by
\begin{itemize}
\item Case $T=0:$ $\gamma^i(\lambda)=\min \{\lambda+\kappa_i,0\}$.
\item Case $T \neq 0, C=0$: $\gamma^i(\lambda)=\frac{T_i W_i}{2 |A_i|} + \min \{\lambda+\kappa_i,0\}$ if $A_i<0$, $\gamma^i(\lambda)=+\infty$ if $A_i \geq 0$.
\item Case $C \neq 0$ and $T \neq 0$:
\begin{equation*}
\gamma^i(\lambda)=\begin{cases}
T_i \, x_{2,i}+\kappa_i+\lambda \; \text{ if } \lambda \leq \lambda_i(x_{2,i}), \\
T_i \, \Sigma^*_{i}(\lambda)+\frac{V_i (\kappa_i+\lambda) (2 A_i \Sigma^*_{i}(\lambda)+W_i)}{C_i^2 (\Sigma^*_{i}(\lambda))^2} \; \text{ if } \lambda_i(x_{2,i}) < \lambda < \lambda_i(x_{e,i}), \\
T_i \, x_{e,i}  \; \text{ if } \lambda_i(x_{e,i}) \leq \lambda.
\end{cases}
\end{equation*}
\end{itemize}
 where in the second case $\Sigma^*_{i}(\lambda)$ is the unique positive root of (\ref{eq: cubic}).
\end{thm}

\begin{proof}
The indexability comes from the fact that the indices $\lambda(\Sigma)$ are verified to be monotonically increasing functions of $\Sigma$. Inverting the relation we obtain $\Sigma_{th}(\lambda)$ as the variance for which we are indifferent between the active and passive actions. As we increase $\lambda$, $\Sigma_{th}(\lambda)$ increases and the passive region (the interval $[0,\Sigma_{th}(\lambda)]$) increases.
\end{proof}

\section{Multidimensional Systems}		\label{section: MIMO case}

Generalizing the computations of the previous section to multidimensional systems requires solving the corresponding optimal control problem in higher dimensions, for which it is not clear that a closed form solution exist. Moreover we have considered in section \ref{section: 1D-dynamics} a particular case of the sensor scheduling problem where all sensors are identical. We now return to the general multidimensional problem and sensors with possibly distinct characteristics, as described in the introduction. 

For the infinite-horizon average cost problem, we show that computing the value of a lower bound similar to the one presented in section \ref{section: 1D-dynamics} reduces to a convex optimization problem involving, at worst, Linear Matrix Inequalities (LMI) whose size grows polynomially with the problem essential parameters. Moreover, one can further decompose the computation of this convex program into $N$ coupled subproblems as in the standard restless bandit case.

\subsection{Performance Bound}

For convenience, let us repeat the deterministic optimal control problem under consideration:
\begin{align}	\label{eq: OCP continuous}
& \min_\pi \; \lim \sup_{T \to \infty} \frac{1}{T} \int_0^T \sum_{i=1}^N \left \{ \trace (T_i \, \Sigma_i(t)) + \sum_{j=1}^{M} \kappa_{ij} \pi_{ij}(t) \right \} dt, \\
& \text{s.t. } \; \dot \Sigma_i(t) = A_i \Sigma_i + \Sigma_i A_i^T+W_i-\Sigma_i \left(  \sum_{j=1}^M \pi_{ij} (t) C_{ij}^T V_{ij}^{-1} C_{ij} \right) \Sigma_i, \; i=1\ldots,N, \label{eq: Riccati diff eq. continuous} \\
& \pi_{ij}(t) \in \{0,1\},  \; \forall t \geq 0,  i=1\ldots,N, \; j=1,\ldots,M, \nonumber \\
& \sum_{i=1}^N \pi_{ij}(t) \leq 1,  \; \forall t \geq 0, \;\; j=1,\ldots,M, \label{eq: constraint continuous} \\
& \sum_{j=1}^M \pi_{ij}(t) \leq 1,  \; \forall t \geq 0, \;\; i=1,\ldots,N, \label{eq: constraint continuous II} \\
& \Sigma_i(0)=\Sigma_{i,0}, \; i=1,\ldots,N. \nonumber
\end{align}
Here we consider the constraints (\ref{eq: pathwise resource constraint}) and (\ref{eq: pathwise resource constraint II}), but any combination of inequality and equality constraints from (\ref{eq: pathwise resource constraint})-(\ref{eq: pathwise resource constraint II bis}) can be used without change in the argument for the derivation of the performance bound.
We define the following quantities:
\begin{equation}	\label{eq: lim sup on policy}
\tilde \pi_{ij}(T) = \frac{1}{T} \int_0^T \pi_{ij}(t) dt, \, \forall \, T \geq 0.
\end{equation}
Since $\pi_{ij}(t) \in \{0,1\}$ we must have $0 \leq \tilde \pi_{ij}(T) \leq 1$. 
Our first goal, inspired by the idea already exploited in the restless bandit problem, is to obtain a lower bound on the cost of the finite-horizon optimal control problem in terms of the numbers $\tilde \pi_{ij}(T)$ instead of the functions $\pi_{ij}(t)$.

It will be easier to work with the information matrices
\[
Q_i(t) = \Sigma^{-1}_i(t).
\]
Note that invertibility of $\Sigma_i(t)$ is guaranteed by our assumptions, as a consequence of \cite[theorem 21.1]{Brockett70_book}. Hence we replace the dynamics (\ref{eq: Riccati diff eq. continuous}) by the equivalent
\begin{equation}	\label{eq: Riccati information}
\dot Q_i = - Q_i A_i - A_i^T Q_i - Q_i W_i Q_i + \sum_{j=1}^M \pi_{ij} (t) C_{ij}^T V_{ij}^{-1} C_{ij}, \quad i=1,\ldots,N.
\end{equation}
Let us also define, for all $T$,
\begin{align*}
\tilde \Sigma_{i}(T) := \frac{1}{T} \int_0^T \Sigma_{i}(t) dt, \quad
\tilde Q_i(T) :=\frac{1}{T} \int_0^T Q_i(t) dt. 
\end{align*}
By linearity of the trace operator, we can rewrite the objective function
\[
\limsup_{T \to \infty} \sum_{i=1}^N \left \{ \trace ( T_i \, \tilde \Sigma_i(T) ) +\sum_{j=1}^M \kappa_{ij} \tilde \pi_{ij}(T) \right \}.
\]
Let $\mathbb{S}^n, \mathbb{S}^n_{+}, \mathbb{S}^n_{++}$ denote the set of symmetric, symmetric positive semidefinite and symmetric positive definite matrices respectively. A function $f: \mathbb{R}^m \to \mathbb{S}^n$ is called \emph{matrix convex} if and only if for all $x,y \in \mathbb{R}^m$ and $\alpha \in [0,1]$, we have
\[
f(\alpha x + (1-\alpha)y) \preceq \alpha f(x) + (1-\alpha) f(y),
\]
where $\preceq$ refers to the usual partial order on $\mathbb{S}^n$, i.e., $A \preceq B$ if and only if $B-A \in \mathbb{S}^n_{+}$. Equivalently, $f$ is matrix convex if the scalar function $x \mapsto z^T f(x) z$ is convex for all vectors $z$. The following lemma will be useful

\begin{lemma}	\label{lem: matrix convexity}
The functions 
\begin{align*}
& \mathbb{S}^n_{++} \to \mathbb{S}^n_{++} & \mathbb{S}^n  \to \mathbb{S}^n \\
& X \mapsto X^{-1} & X \mapsto XWX
\end{align*}
for $W \in \mathbb{S}^n_{+}$, are matrix convex.
\end{lemma}

\begin{proof}
See \cite[p.76, p.110]{boyd_cvxBook06}.
\end{proof}

A consequence of this lemma is that Jensen's inequality is valid for these functions. We use it first as follows
\[
\forall T, \;\; \left( \frac{1}{T} \int_0^T \Sigma_i(t) dt \right)^{-1} \preceq \frac{1}{T} \int_0^T Q_i(t) dt = \tilde Q_i(T), 
\]
hence
\begin{align*}
\forall T, \;\; \tilde \Sigma_i(T) \succeq (\tilde Q_i(T))^{-1}.
\end{align*}
and so
\[
\trace (T_i \, \tilde \Sigma_i(T)) \geq \trace (T_i \, (\tilde Q_i(T))^{-1}).
\]
Next, integrating (\ref{eq: Riccati information}) and letting $Q_{i,0}=\Sigma_{i,0}^{-1}$, we have
\begin{align}
\frac{1}{T} (Q_i(T)-Q_{i,0}) &= - \tilde Q_i(T) A_i - A_i^T \tilde Q_i(T) - \frac{1}{T} \int_0^T Q_i(t) W_i Q_i(t) dt + \sum_{j=1}^M \left( \frac{1}{T} \int_0^T \pi_{ij}(t) \right) C_{ij}^T V_{ij}^{-1} C_{ij} \nonumber \\
\frac{1}{T} (Q_i(T)-Q_{i,0}) &= - \tilde Q_i(T) A_i - A_i^T \tilde Q_i(T) - \frac{1}{T} \int_0^T Q_i(t) W_i Q_i(t) dt + \sum_{j=1}^M \tilde \pi_{ij}(T) C_{ij}^T V_{ij}^{-1} C_{ij} \nonumber .
\end{align}

Using Jensen's inequality and lemma \ref{lem: matrix convexity} again, we have
\[
\frac{1}{T} \int_0^T Q_i(t) W_i Q_i(t) \succeq \tilde Q_i(T) W_i \tilde Q_i(T),
\]
and so we obtain
\begin{align}
\frac{1}{T} (Q_i(T)-Q_{i,0}) & \preceq - \tilde Q_i(T) A_i - A_i^T \tilde Q_i(T) - \tilde Q_i(T) W_i \tilde Q_i(T) + \sum_{j=1}^M \tilde \pi_{ij}(T) C_{ij}^T V_{ij}^{-1} C_{ij}. \label{eq: Riccati integrated}
\end{align}
Last, since $Q_i(T) \succeq 0$, this implies, for all $T$,
\begin{equation}	\label{eq: Riccati integrated 2}
\tilde Q_i(T) A_i + A_i^T \tilde Q_i(T) + \tilde Q_i(T) W_i \tilde Q_i(T) - \sum_{j=1}^M \tilde \pi_{ij}(T) C_{ij}^T V_{ij}^{-1} C_{ij} \preceq \frac{Q_{i,0}}{T}.
\end{equation}

So we see that for a fixed policy $\pi$ and any time $T$, the quantity
\begin{equation}	\label{eq: objective finite T}
\sum_{i=1}^N \left \{ \trace ( T_i \, \tilde \Sigma_i(T) ) +\sum_{j=1}^M \kappa_{ij} \tilde \pi_{ij}(T) \right \}
\end{equation}
is lower bounded by the quantity
\[
\sum_{i=1}^N \left \{ \trace ( T_i \, ( \tilde Q_i(T))^{-1} ) +\sum_{j=1}^M \kappa_{ij} \tilde \pi_{ij}(T) \right \},
\]
where the matrices $\tilde Q_i(T)$ and the number $\tilde \pi_{ij}(T)$ are subject to the constraints $(\ref{eq: Riccati integrated 2})$ as well as
\[
0 \leq \tilde \pi_{ij}(T) \leq 1, \;\; \sum_{i=1}^N \tilde \pi_{ij}(T) \leq 1,  \; j=1,\ldots,M,  \;\; \sum_{j=1}^M \tilde \pi_{ij}(T) \leq 1,  \; i=1,\ldots,N. \nonumber
\]
Hence for any $T$, the quantity $Z^*(T)$ defined below is a lower bound on the value of $(\ref{eq: objective finite T})$ for any choice of policy $\pi$
\begin{align}	
& Z^*(T) = \min_{ \tilde Q_i, p_{ij}} \quad \sum_{i=1}^N \left \{ \trace (T_i \, \tilde Q^{-1}_i) + \sum_{j=1}^{M} \kappa_{ij} p_{ij} \right \}, \label{eq: optimization problem, fixed T} \\
& \text{s.t. } \tilde Q_i A_i + A_i^T \tilde Q_i + \tilde Q_i W_i \tilde Q_i - \sum_{j=1}^M p_{ij} C_{ij}^T V_{ij}^{-1} C_{ij} \preceq \frac{Q_{i,0}}{T}, \; i=1\ldots,N, \label{eq: Riccati inequality with RHS} \\
& \tilde Q_i \succ 0, \; i=1\ldots,N, \nonumber \\
& 0 \leq p_{ij} \leq 1,  \; i=1 \ldots,N, \; j=1,\ldots,M, \nonumber \\
& \sum_{i=1}^N p_{ij} \leq 1,  \; j=1,\ldots,M,  \nonumber \\
& \sum_{j=1}^M p_{ij} \leq 1,  \; i=1,\ldots,N. \nonumber
\end{align}
Consider now the following program, where the right-hand side of (\ref{eq: Riccati inequality with RHS}) has been replaced by $0$:
\begin{align}	
& Z^* = \min_{ \tilde Q_i, p_{ij}} \quad \sum_{i=1}^N \left \{ \trace (T_i \, \tilde Q^{-1}_i) + \sum_{j=1}^{M} \kappa_{ij} p_{ij} \right \}, \label{eq: optimal value} \\
& \text{s.t. } \tilde Q_i A_i + A_i^T \tilde Q_i + \tilde Q_i W_i \tilde Q_i - \sum_{j=1}^M p_{ij} C_{ij}^T V_{ij}^{-1} C_{ij} \preceq 0, \; i=1\ldots,N, \label{eq: Riccati inequality} \\
& \tilde Q_i \succ 0, \; i=1\ldots,N, \nonumber \\
& 0 \leq p_{ij} \leq 1,  \; i=1 \ldots,N, \; j=1,\ldots,M, \nonumber \\
& \sum_{i=1}^N p_{ij} \leq 1,  \; j=1,\ldots,M,  \nonumber \\
& \sum_{j=1}^M p_{ij} \leq 1,  \; i=1,\ldots,N. \nonumber
\end{align}
Defining $\delta:=1/T$, and rewriting with a slight abuse of notation $Z^*(\delta)$ instead of $Z^*(T)$ for $\delta$ positive, we also define $Z^*(0)=Z^*$, where $Z^*$ is given by (\ref{eq: optimal value}). Note that $Z^*(0)$ is finite, since we can find a feasible solution as follows. 
For each $i$, we choose a set of indices $J_i=\{j_1, \ldots, j_{n_i}\} \subset \{1,\ldots,M\}$ such that $(A_i,\tilde C_{i})$ is observable, as in assumption \ref{assumption: observability}. Once a set $J_i$ has been chosen for each $i$, we form the matrix $\hat P$ with elements $\hat p_{ij}=1_{\{j \in J_i\}}$. Finally, we form a matrix $P$ with elements $p_{ij}$ satisfying the constraints and nonzero exactly where the $\hat p_{ij}$ are nonzero. 
Such a matrix is easy to find if we consider the inequality constraints (\ref{eq: pathwise resource constraint}) and (\ref{eq: pathwise resource constraint II}). If equality constraints are involved instead, such a matrix $P$ exists as a consequence of Birkhoff theorem \cite{birkhoff46_3obs}, see theorem \ref{thm: Mirsky thm}. Now we consider the quadratic inequality (\ref{eq: Riccati inequality}) for some value of $i$. From the detectability assumption \ref{assumption: observability} and the choice of $p_{ij}$, we deduce that the pair $(A_i, \hat C_i)$, with
\begin{equation}	\label{eq: composite C matrix}
\hat C_i =  \left [ \begin{array}{ccc} \sqrt{p_{i 1}} \, C^T_{i 1} V_{i 1}^{-1/2} & \cdots & \sqrt{p_{i M}} \, C^T_{i M} V_{i M}^{-1/2} \end{array} \right ]^T
\end{equation}
is detectable. Also note that
\[
\hat C^T_i \hat C_i = \sum_{j=1}^{M} p_{ij} \, C_{ij}^T V_{ij}^{-1} C_{ij}. 
\]
Together with the controllability assumption \ref{assumption: controllability}, we then know that (\ref{eq: Riccati inequality}) has a positive definite solution $\tilde Q_i$ \cite[theorem 2.4.25]{Aboukandil03_mreBook}. Hence $Z^*(0)$ is finite. 

We can also define $Z^*(\delta)$ for $\delta<0$, by changing the right-hand side of  (\ref{eq: Riccati inequality with RHS}) into $\delta Q_{i,0}=-|\delta| Q_{i,0}$. We have that $Z^*(\delta)$ is finite for $\delta<0$ small enough. Indeed, passing the term $\delta Q_{i,0}$ on the left hand side, this can then be seen as a perturbation of the matrix $\hat C_i$ above, and for $\delta$ small enough, detectability, which is an open condition, is preserved.
Now we will see below that (\ref{eq: optimization problem, fixed T}), (\ref{eq: optimal value}) are convex programs. It is then a standard result of perturbation analysis (see e.g. \cite[p. 250]{boyd_cvxBook06}) that $Z^*(\delta) $ is a convex function of $\delta$, hence continuous on the interior of its domain, in particular continuous at $\delta=0$. So
\[
\limsup_{T \to \infty} Z^*(T) = \lim_{T \to \infty} Z^*(T) = Z^*.
\]
Finally, for any policy $\pi$, we obtain the following lower bound on the achievable cost
\[
\limsup_{T \to \infty} \frac{1}{T} \int_0^T \sum_{i=1}^N \left \{ \trace ( T_i \, \Sigma_i(t) ) +\sum_{j=1}^M \kappa_{ij} \pi_{ij}(t) \right \} dt \geq \lim_{T \to \infty} Z^*(T) = Z^*.
\]
We now show how to compute $Z^*$ by solving a convex program involving linear matrix inequalities. For each $i$, introduce a new (slack) matrix variable $R_i$. Since $Q_i \succ 0, \, R_i \succ Q_i^{-1}$ is equivalent, by taking the Schur complement, to
\[
\left[
\begin{array}{cc}
R_i & I \\
I & Q_i
\end{array}
\right]
\succ 0,
\]
and the Riccati inequality (\ref{eq: Riccati inequality}) can be rewritten
\[
\left[
\begin{array}{cc}
\tilde Q_i A_i + A_i^T \tilde Q_i - \sum_{j=1}^M p_{ij} C_{ij}^T V_{ij}^{-1} C_{ij} & \tilde Q_i W_i^{1/2} \\
W_i^{1/2} \tilde Q_i & -I
\end{array}
\right] \preceq 0.
\]
We finally obtain, dropping the tildes from the notation $\tilde Q_i$, the semidefinite program
\begin{align}	\label{eq: LMI relaxation}
& Z^* = \min_{ R_i, Q_i, p_{ij}} \quad \sum_{i=1}^N \left \{ \trace (T_i \, R_i) + \sum_{j=1}^{M} \kappa_{ij} p_{ij} \right \}, \\
& \text{s.t. } \left[
\begin{array}{cc}
R_i & I \\
I & Q_i
\end{array}
\right]
\succ 0, \; i=1\ldots,N,  \nonumber \\
& \left[
\begin{array}{cc}
Q_i A_i + A_i^T Q_i - \sum_{j=1}^M p_{ij} C_{ij}^T V_{ij}^{-1} C_{ij} & Q_i W_i^{1/2} \\
W_i Q_i^{1/2} & -I
\end{array}
\right] \preceq 0, \; i=1\ldots,N, \nonumber \\ 
& 0 \leq p_{ij} \leq 1,  \; i=1 \ldots,N, \; j=1,\ldots,M, \nonumber \\
& \sum_{i=1}^N p_{ij} \leq 1,  \; j=1,\ldots,M, \label{eq: coupling constraint LMI} \\ 
& \sum_{j=1}^M p_{ij} \leq 1,  \; i=1,\ldots,N. \nonumber 
\end{align}
Hence solving the program (\ref{eq: LMI relaxation}) provides a lower bound on the achievable cost for the original optimal control problem.

\subsection{Problem Decomposition}

It is well-know that efficient methods exist to solve (\ref{eq: LMI relaxation}) in polynomial time, which implies a computation time polynomial in the number of variables of the original problem. Still, as the number of targets increases, the large LMI (\ref{eq: LMI relaxation}) becomes difficult to solve. Note however that it can be decomposed into $N$ small coupled LMIs, following the standard dual decomposition approach already used for the restless bandit problem. This decomposition is very useful to solve large scale programs with a large number of systems. For completeness, we present the argument in more details below.

We first note that (\ref{eq: coupling constraint LMI}) is the only constraint which links the $N$ subproblems together. So we form the Lagrangian
\[
L(R,Q,p;\lambda) = \sum_{i=1}^N \left \{ \trace (T_i \, R_i) + \sum_{j=1}^{M} ( \kappa_{ij} + \lambda_j ) p_{ij} \right \} - \sum_{j=1}^M \lambda_j,
\]
where $\lambda \in \R_+^M$ is a vector of Lagrange multipliers. We would take $\lambda \in \R^M$ if we had the constraint (\ref{eq: pathwise resource constraint bis}) instead of (\ref{eq: pathwise resource constraint}). Now the dual function is
\begin{equation}	\label{eq: dual function value LMI}
G(\lambda)=\sum_{i=1}^N G_i(\lambda)-\sum_{j=1}^M \lambda_j,
\end{equation}
\begin{align}
\text{with } \quad G_i(\lambda) & = \min_{ R_i, Q_i, \{p_{ij}\}_{1 \leq j \leq M}} \quad \trace (T_i \, R_i) + \sum_{j=1}^{M} (\kappa_{ij}+\lambda_j) p_{ij},   \label{eq: local LMI} \\
& \text{s.t. } \left[
\begin{array}{cc}
R_i & I \\
I & Q_i
\end{array}
\right]
\succ 0,  \nonumber \\
& \left[
\begin{array}{cc}
Q_i A_i + A_i^T Q_i - \sum_{j=1}^M p_{ij} C_{ij}^T V_{ij}^{-1} C_{ij} & Q_i W_i^{1/2} \\
W_i^{1/2} Q_i & -I
\end{array}
\right] \preceq 0, \nonumber \\ 
& 0 \leq p_{ij} \leq 1, \; j=1,\ldots,M, \nonumber \\
& \sum_{j=1}^M p_{ij} \leq 1. \nonumber 
\end{align}

The optimization algorithm proceeds then as follows \cite[chap. 11]{bonnans06_book}. We choose an initial value $\lambda^1 \geq 0$ and set $k=1$.
\begin{enumerate}
\item For $i=1,\ldots,N$, compute  $R^k_i, Q^k_i,\{ p^k_{ij} \}_{1 \leq j \leq M}$ optimal solution of (\ref{eq: local LMI}), and the value $G_i(\lambda^k)$.	\label{subgradient: first step}
\item The value of the dual function at $\lambda^k$ is given by (\ref{eq: dual function value LMI}). A supergradient of $G(\lambda^k)$ at $\lambda^k$ is given by
\[
\left [ \sum_{i=1}^N p^k_{i1} - 1, \ldots,  \sum_{i=1}^N p^k_{iM} - 1 \right ]. 
\]
\item Compute $\lambda^{k+1}$ in order to maximize $G(\lambda)$. We can do this by using a supergradient algorithm, or any preferred nonsmooth optimization algorithm. Let k:=k+1 and go to step \ref{subgradient: first step}, or stop if convergence is satisfying.
\end{enumerate}

Because the initial program (\ref{eq: LMI relaxation}) is convex, we know that the optimal value of the dual optimization problem is equal to the optimal value of the primal. Moreover, the optimal variables of the primal are obtained at step \ref{subgradient: first step} of the algorithm above once convergence has been reached.

\subsection{Open-loop Periodic Policies Achieving the Performance Bound}		\label{section: BvN policies}

\subsubsection{Definition of the Policies}

In this section we describe a sequence of open-loop policies that can approach arbitrarily closely the lower bound computed by (\ref{eq: LMI relaxation}), thus proving that this bound is tight. These policies are periodic switching strategies using a schedule obtained from the optimal parameters $p_{ij}$. Assuming no switching times or costs, their performance approaches the bound as the length of the switching cycle decreases toward $0$.

Let $P=[p_{ij}]_{1 \leq i \leq N, 1 \leq j \leq M}$ be the matrix of optimal parameters obtained in the solution of (\ref{eq: LMI relaxation}). We assume here that constraints (\ref{eq: pathwise resource constraint}) and (\ref{eq: pathwise resource constraint II}) were enforced, which is the most general case for the discussion in this section.  Hence $P$ verifies
\begin{align*}
& 0 \leq p_{ij} \leq 1, \; i=1,\ldots,N, \; j=1,\ldots,M, \\
& \sum_{i=1}^N p_{ij} \leq 1, \; j=1,\ldots,M, \;\; \text{ and } \sum_{j=1}^M p_{ij} \leq 1, \; i=1,\ldots,N.
\end{align*}
A \emph{doubly substochastic matrix} of dimension $n$ is an $n \times n$ matrix $A=[a_{ij}]_{1 \leq i,j \leq n}$ which satisfies
\begin{align*}
& 0 \leq a_{ij} \leq 1, \; i,j=1,\ldots,n, \\
& \sum_{i=1}^n a_{ij} \leq 1, \; j=1,\ldots,n, \;\; \text{ and } \sum_{j=1}^n a_{ij} \leq 1, \; i=1,\ldots,n.
\end{align*}
If $M=N$, $P$ is therefore a doubly substochastic matrix. Else if $M<N$ (resp. $N<M$) we can add $N-M$ columns of zeros (resp. $M-N$ rows of zeros) to $P$ to obtain a doubly substochastic matrix. In any case, we call the resulting doubly substochastic matrix $\tilde P=[\tilde p_{ij}]$. If rows have been added, this is equivalent to the initial problem with additional ``dummy systems''. If columns are added, these correspond to using ``dummy sensors''. Dummy systems (i.e., for $i>N$) are not included in the objective function (the corresponding $T_i$ is $0$), and a dummy sensor (i.e., for $j>M$) is associated formally to the measurement noise covariance matrix $V_{ij}^{-1}=0$ for all $i$, in effect producing no measurement. In the following we assume that $\tilde P$ is an $N \times N$ doubly substochastic matrix, but the discussion in the $M \times M$ case is identical.
Doubly substochastic matrices have been intensively studied, and the material used in the following can be found in the book of Marshall and Olkin \cite{marshall79_book}. In particular, we have the following corollary of a classical theorem of Birkhoff \cite{birkhoff46_3obs}, which says that a doubly stochastic matrix is a convex combination of permutation matrices.
\begin{theorem}[\cite{Mirsky59_matrices}]			\label{thm: Mirsky thm}
The set of $N \times N$ doubly substochastic matrices is the convex hull of the set $\mathcal P_0$ of $N \times N$ matrices which have a most one unit in each row and each column, and all other entries are zero.
\end{theorem}

Hence for the doubly substochastic matrix $\tilde P$, there exists a set of positive numbers $\phi_k$ and matrices $P_k \in \mathcal P_0$ such that 
\begin{align}	\label{eq: Mirsky's decomposition}
\tilde P = \sum_{k=1}^K \phi_k P_k, \; \text{ with }
\sum_{k=1}^K \phi_k = 1, \; \text{ for some integer } K.
\end{align}
One way of computing this decomposition is to first extend $\tilde P$ to the $2N \times 2N$ doubly stochastic matrix
\[
\hat P = \left[
\begin{array}{cc}
\tilde P & I-D_r \\
I-D_c & \tilde P^T
\end{array}
\right],
\]
where $r_1,\ldots,r_N$ and $c_1,\ldots,c_N$ are the row sums and column sums of $\tilde P$, and $D_r=\text{diag}(r_1,\ldots,r_N)$, $D_c=\text{diag}(c_1,\ldots,c_N)$. Then there is an algorithm that runs in time $O(N^{4.5})$ \cite{marshall79_book, chang00_birkhoffSwitch} and provides the decomposition 
\[
\hat P = \sum_{k=1}^K \phi_k \hat P_k,
\]
with $K=(2N-1)^2+1$ and where the $\hat P_k$'s are permutation matrices of size $2N \times 2N$. The decomposition (\ref{eq: Mirsky's decomposition}) is finally obtained by deleting the last $N$ rows and columns of $\hat P_k$ to obtain the matrices $P_k$, $k=1,\ldots,K$.

Note that any matrix $A=[a_{ij}]_{i,j} \in \mathcal P_0$ represents a valid sensor/system assignment (for the system with additional dummy systems or sensors), where sensor $j$ is measuring system $i$ if and only if $a_{ij}=1$. With the decomposition (\ref{eq: Mirsky's decomposition}), we now consider a family of periodic switching policies parametrized by a positive number $\epsilon$ representing a time interval over which the switching schedule is executed completely. For a given value of $\epsilon$, the policy is defined as follows:
\begin{enumerate}
\item \label{switching step 1} At time $t=l \epsilon, l \in \mathbb{N}$, associate sensor $j$ to system $i$ as specified by the matrix $P_1$ of the representation (\ref{eq: Mirsky's decomposition}). 
 Run the corresponding continuous-time Kalman filters, keeping this sensor/system association for a duration $\phi_1 \epsilon$.
\item At time $t=(l+\phi_1)\epsilon$, switch to the assignment specified by $P_2$. Run the corresponding continuous time Kalman filters until $t=(l+\phi_1+\phi_2)\epsilon$.
\item Repeat the switching procedure, switching to matrix $P_{i+1}$ at time $t=l+\phi_1+\cdots+\phi_i$, for $i=1,\ldots,K-1$.
\item At time $t=(l+\phi_1+\cdots+\phi_K)\epsilon = (l+1)\epsilon$, start the switching sequence again at step \ref{switching step 1} with $P_1$ and repeat the steps above.
\end{enumerate}
It is easy to see that the matrices $P_i, \, i=1,\ldots,K$ never specify that a ``dummy sensor'' should execute a measurement or that a ``dummy system'' should be measured, since from the decomposition (\ref{eq: Mirsky's decomposition}) this would correspond to nonzero entries in the columns or rows added to $P$ to form $\tilde P$.


\subsubsection{Performance of the Periodic Switching Policies}       \label{section: switching policies performance - KF}

Let us fix $\epsilon>0$ in the definition of the switching policy, and consider now, for this policy, the evolution of the covariance matrix $\Sigma^\epsilon_i(t)$ for the estimation error on the state of system $i$. The superscript indicates the dependence on the period $\epsilon$ of the policy. First we have
\begin{lemma}	\label{lem: convergence periodic}
For all $i \in \{1,\ldots,N\}$, the estimation error covariance $\Sigma^\epsilon_i(t)$ converges as $t \to \infty$ to a periodic function $\bar \Sigma^\epsilon_i(t)$ of period $\epsilon$.
\end{lemma}

\begin{proof}
Fix $i \in \{1,\ldots,N\}$. 
Let $\sigma_i(t) \in \{0,1,\ldots,N\}$ be the function specifying which sensor is observing system $i$ at time $t$ under the switching policy. By convention $\sigma_i(t)=0$ means that no sensor is scheduled to observe system $i$, and $\sigma_i(t)=j$ means that sensor $j$ measures system $i$. Note from the remark following the description of the switching policies that in fact we have $\sigma_i(t) \in \{0, \ldots, M\}$, i.e., the policy never schedules measurements by dummy sensors. Similarly, if instead we were considering the situation $M>N$ and $\tilde P$ an $M \times M$ matrix, then we would have $\sigma_i(t)=0$ for $i \in \{N+1,\ldots,M\}$ and all $t$.
Note also that $\sigma_i(t)$ is a piecewise constant, $\epsilon$-periodic function. The switching times of $\sigma_i(t)$ are $t=(l+\phi_1+\cdots+\phi_{k-1})\epsilon$, for $k=1,\ldots,K$ and $l \in \mathbb{N}$.

The covariance matrix $\Sigma^\epsilon_i(t)$ obeys the following periodic Riccati differential equation (PRE):
\begin{align}
\dot \Sigma^\epsilon_i(t) &= A_i \Sigma^\epsilon_i(t) + \Sigma^\epsilon_i(t) A_i^T + W_i - \Sigma^\epsilon_i(t) (C_i^\epsilon(t))^T  C_i^\epsilon(t) \Sigma^\epsilon_i(t) \label{eq: PRE} \\
\Sigma^\epsilon_i(t) &= \Sigma_{i,0}, \nonumber
\end{align}
where $C_i^\epsilon(t) := V_{i \sigma_i(t)}^{-1/2} C_{i \sigma_i(t)}$ is a piecewise constant, $\epsilon$-periodic matrix valued function, and we use the convention $V_{i j}^{-1} = C_{i j}=0$ when $j=0$. 
We now show that $(A_i, C_i^\epsilon(\cdot))$ is detectable. Let $j_1, \ldots, j_K$ be the successive values taken by the function $\sigma_i(t)$ over the period $\epsilon$. From the definition of detectability for linear periodic systems and it modal characterization \cite[p.130]{Bittanti91_PRE}, we immediately deduce that the pair $(A_i,C_i^\epsilon(\cdot))$ is not detectable if and only if there exists an eigenpair $(\lambda,x)$ for $A_i$, with $\text{Re}(\lambda) \geq 0$, $x \neq 0$, such that
\begin{flalign}
&& A_i x =& \lambda x, \text{ and } C^\epsilon_i(t) e^{A_it}x=e^{\lambda t} C^\epsilon_i(t) x = 0, \, \forall t \in [0,\epsilon], && \nonumber \\
\textrm{hence } && C_{i j_1}x =& \ldots=C_{i j_K}x=0. && \label{eq: modal test}
\end{flalign}
Let us denote by $p_{k,ij}$ the $(i,j)^{\text{th}}$ element of the matrix $P_k$ in the decomposition (\ref{eq: Mirsky's decomposition}). We have $p_{k,ij}=1_{ \{j=j_k \} }$ with the above definition of $j_k$, including the case $j_k=0$ (no measurement), which gives $p_{k,ij}=0$ for all $j \in \{1,\ldots,N\}.$
Then we can write
\begin{align}
\sum_{k=1}^K \phi_k C^T_{i j_k} V^{-1}_{i j_k} C_{i j_k} &= \sum_{k=1}^K \phi_k \left( \sum_{j=1}^N p_{k,ij} C^T_{ij} V^{-1}_{ij} C_{ij} \right)  \nonumber \\
&=\sum_{j=1}^N \left( \sum_{k=1}^K \phi_k p_{k,ij} \right) C^T_{ij} V^{-1}_{ij} C_{ij} \nonumber \\
&=\sum_{j=1}^N \tilde p_{ij} C^T_{ij} V^{-1}_{ij} C_{ij}  \nonumber \\
&=\sum_{j=1}^M p_{ij} C^T_{ij} V^{-1}_{ij} C_{ij} = \hat C_i^T \hat C_i,	\label{eq: rewriting as frequencies}
\end{align}
where the next-to-last equality uses the fact that $\tilde p_{ij} = p_{ij}$ for $j \leq M$ and $\tilde p_{ij} = 0$ for $j \geq M+1$,
and $\hat C_i$ was defined in (\ref{eq: composite C matrix}). Note that we now consider this definition of $\hat C_i$ for the optimal parameters $p_{ij}$ provided by the solution of (\ref{eq: LMI relaxation}). Then (\ref{eq: modal test}) and (\ref{eq: rewriting as frequencies}) imply $\| \hat C_i x \|_2 = 0$, so $\hat C_i x = 0$, i.e., $(A, \hat C_i)$ is not detectable.
But the parameters $p_{ij}$ being optimal for the program (\ref{eq: LMI relaxation}), this would imply that this program is not feasible  \cite[p.68]{Aboukandil03_mreBook}, a contradiction with our discussion following (\ref{eq: composite C matrix}).
So $(A_i,C_i^\epsilon(\cdot))$ must be detectable, and together with our assumption \ref{assumption: pos-def initial condition}, this implies by the result of \cite[p. 95]{denicolao92_convergence} that 
\[
\lim_{t \to \infty} \left ( \Sigma_i^\epsilon(t) - \bar \Sigma^\epsilon_i(t) \right) = 0,
\]
where $\bar \Sigma^\epsilon_i(t)$ is the strong solution of the PRE, which is $\epsilon$-periodic.

\end{proof}

Next, denote by $\tilde \Sigma_i(t)$ the solution to the following Riccati differential equation (RDE):
\begin{align} \label{eq: averaged RDE}
\dot \Sigma_i = A_i \Sigma_i + \Sigma_i A_i^T + W_i - \Sigma_i \left( \sum_{j=1}^M p_{ij} C_{ij}^T V_{ij}^{-1} C_{ij} \right) \Sigma_i \,, \quad \Sigma_i(0)=\Sigma_{i,0}.
\end{align}
Assumptions \ref{assumption: observability} and \ref{assumption: controllability}, together with our discussion of the implied detectability of the pair $(A,\hat C_i)$ (see (\ref{eq: rewriting as frequencies})), guarantee that $\tilde \Sigma_i(t)$ converges to a positive definite limit denoted $\Sigma^*_i$. Moreover, $\Sigma_i^*$ is stabilizing and is the unique positive definite solution to the algebraic Riccati equation (ARE):
\begin{equation}	\label{eq: averaged ARE - base}
A_i \Sigma_i + \Sigma_i A_i^T + W_i - \Sigma_i \left( \sum_{j=1}^M p_{ij} C_{ij}^T V_{ij}^{-1} C_{ij} \right) \Sigma_i = 0.
\end{equation}


The next lemma says that the periodic function $\bar \Sigma^\epsilon_i(t)$ oscillates in a neighborhood of $\Sigma_i^*$.

\begin{lemma}	\label{lem: approximation of average}
For all $t \in \R_+$, we have $\bar \Sigma^\epsilon_i(t)-\Sigma_i^* = O(\epsilon)$ as $\epsilon \to 0$. 
\end{lemma}

\begin{proof}
The function $t \to \bar \Sigma^\epsilon_i(t)$ of lemma \ref{lem: convergence periodic} is the strong periodic solution of the PRE (\ref{eq: PRE}). It is $\epsilon$-periodic. From Radon's lemma  \cite[p.90]{Aboukandil03_mreBook}, which gives a representation of the solution to a Riccati differential equation as the ratio of solutions to a linear ODE, we also know that $\bar \Sigma_i^\epsilon$ is $C^\infty$ on each interval where $\sigma_i(t)$ is constant, where $\sigma_i(t)$ is the switching signal defined in the proof of lemma \ref{lem: convergence periodic}.

Let $\hat \Sigma^\epsilon_i$ be the average of $t \to \bar \Sigma^\epsilon_i(t)$:
\[
\hat \Sigma^\epsilon_i = \frac{1}{\epsilon} \int_0^\epsilon \bar \Sigma^\epsilon_i(t) dt.
\]
From the preceding remarks, it is easy to deduce that for all $t$, we have $\bar \Sigma^\epsilon_i(t) - \hat \Sigma^\epsilon_i=O(\epsilon).$ Now, averaging the PRE (\ref{eq: PRE}) over the interval $[0,\epsilon]$, we obtain
\[
A_i \hat \Sigma^\epsilon_i + \hat \Sigma^\epsilon_i A^T + W_i - \frac{1}{\epsilon} \int_0^\epsilon  \bar \Sigma^\epsilon_i(t) (C^\epsilon_i(t))^T C^\epsilon_i(t) \bar \Sigma^\epsilon_i(t) dt = \frac{1}{\epsilon}(\bar \Sigma^\epsilon_i(\epsilon) - \bar \Sigma^\epsilon_i(0)) = 0,
\]
where $C^\epsilon_i(t)$ was defined below display (\ref{eq: PRE}).
Expanding this equation in powers of $\epsilon$, we get
\[
A_i \hat \Sigma^\epsilon_i + \hat \Sigma^\epsilon_i A^T + W_i - \hat \Sigma^\epsilon_i  \left( \frac{1}{\epsilon} \int_0^\epsilon (C^\epsilon_i(t))^T C^\epsilon_i(t) dt \right) \hat \Sigma^\epsilon_i + R(\epsilon) = 0,
\]
where $R(\epsilon)=O(\epsilon)$. Let $j_k:=\sigma_i(t)$ for $t \in [(l+\phi_1+\ldots+\phi_{k-1})\epsilon,(l+\phi_1+\ldots+\phi_{k})\epsilon]$. We can then rewrite, using (\ref{eq: rewriting as frequencies}),
\[
\frac{1}{\epsilon} \int_0^\epsilon (C^\epsilon_i(t))^T C^\epsilon_i(t) dt = \sum_{k=1}^K \phi_k C_{i j_k}^T V_{i j_k}^{-1} C_{i j_k}
= \sum_{j=1}^M p_{ij} C^T_{ij} V^{-1}_{ij} C_{ij}.
\]
So we obtain
\[
A_i \hat \Sigma^\epsilon_i + \hat \Sigma^\epsilon_i A^T + (W_i + R(\epsilon)) - \hat \Sigma^\epsilon_i  \left( \sum_{j=1}^M p_{ij} C^T_{ij} V^{-1}_{ij} C_{ij} \right) \hat \Sigma^\epsilon_i = 0.
\] 
Note moreover that for $\epsilon$ sufficiently small, $\hat \Sigma^\epsilon_i$ is the unique positive definite stabilizing solution of this ARE, using the fact that controlability of $(A_i,W_i^{1/2})$ is an open condition. Now comparing this ARE to the ARE (\ref{eq: averaged ARE - base}), and since the stabilizing solution of an ARE is a real analytic function of the parameters \cite{Delchamps84_TACafc}, we deduce that $\hat \Sigma^\epsilon_i-\Sigma^*_i = O(\epsilon)$, and the lemma.

\end{proof}


\begin{theorem}	
Let $Z^\epsilon$ denote the performance of the periodic switching policy with period $\epsilon$. Then $Z^\epsilon - Z^*=O(\epsilon)$ as $\epsilon \to 0$, where $Z^*$ is the performance bound (\ref{eq: LMI relaxation}). Hence the switching policy approaches the lower bound arbitrarily closely as the period tends to $0$.
\end{theorem}

\begin{proof}
We have
\begin{align*}
Z^\epsilon & = \limsup_{T \to \infty} \frac{1}{T} \int_{0}^T \sum_{i=1}^N \left( \trace (T_i \, \Sigma^\epsilon_i(t)) + \sum_{j=1}^M \kappa_{ij} \pi^\epsilon_{ij}(t) \right) dt,
\end{align*}
where $\pi^\epsilon$ is the sensor/system assignment of the switching policy. First by using a transformation similar to (\ref{eq: rewriting as frequencies}) and using the convention $\kappa_{ij}=0$ for $j \in \{0\} \cup \{M+1,\ldots,N\}$ (no measurement or measurement by a dummy sensor), we have for system $i$
\begin{align}
\frac{1}{T} \int_0^T \sum_{j=1}^M \kappa_{ij} \pi_{ij}(t) dt &= \frac{1}{T} \sum_{n=0}^{ \lfloor \frac{T}{\epsilon}  \rfloor -1} \int_{n \epsilon}^{(n+1)\epsilon} \sum_{j=1}^M \kappa_{ij} \pi_{ij}(t) dt + \frac{1}{T} \int_{\lfloor \frac{T}{\epsilon} \rfloor \epsilon}^T \sum_{j=1}^M \kappa_{ij} \pi_{ij}(t) dt    \nonumber \\
&= \frac{1}{T} \left \lfloor \frac{T}{\epsilon} \right \rfloor \sum_{k=1}^K \kappa_{i j_k} (\phi_{k} \epsilon) + \frac{1}{T} \int_{\lfloor \frac{T}{\epsilon} \rfloor \epsilon}^T \sum_{j=1}^M \kappa_{ij} \pi_{ij}(t) dt    \nonumber \\
& = \left \lfloor \frac{T}{\epsilon} \right \rfloor \frac{\epsilon}{T} \sum_{j=1}^M  \kappa_{ij} p_{ij} +\frac{1}{T}  \int_{\lfloor \frac{T}{\epsilon} \rfloor \epsilon}^T \sum_{j=1}^M \kappa_{ij} \pi_{ij}(t) dt,	\label{eq: measurement cost periodic}
\end{align}
where the $j_k$'s were defined in the proof of lemma \ref{lem: convergence periodic}. Hence
\[
 \limsup_{T \to \infty} \frac{1}{T} \int_{0}^T \sum_{i=1}^N \sum_{j=1}^M \kappa_{ij} \pi^\epsilon_{ij}(t) dt \leq \sum_{i=1}^N \sum_{j=1}^M \kappa_{ij} p_{ij}.
\]

Next, from lemma \ref{lem: convergence periodic} and \ref{lem: approximation of average}, it follows readily that $\limsup_{t \to \infty} \Sigma^\epsilon_i(t) - \Sigma_i^* = O(\epsilon)$. It is well known \cite{willems91_ARE} that under our assumptions $\Sigma_i^*$ is the minimal (for the partial order on $\mathbb{S}^n_+$) positive definite solution of the quadratic matrix inequality
\[
A_i \Sigma_i + \Sigma_i A_i^T + W_i - \Sigma_i \left( \sum_{j=1}^M p_{ij} C^T_{ij} V^{-1}_{ij} C_{ij}\right) \Sigma_i \preceq 0.
\]
Hence for the $p_{ij}$ obtained from the computation of the lower bound (\ref{eq: LMI relaxation}), these matrices $\Sigma_i^*$ minimize
\begin{align}
& \sum_{i=1}^N \trace \; (T_i \, \Sigma_i) \nonumber \\
\text{s.t. } \quad & A_i \Sigma_i + \Sigma_i A_i^T + W_i - \Sigma_i \left( \sum_{j=1}^M p_{ij} C^T_{ij} V^{-1}_{ij} C_{ij}\right) \Sigma_i \preceq 0, \quad i=1,\ldots,N. \label{eq: inequalities Riccati sigma}
\end{align}
Changing variable to $Q_i= \Sigma_i^{-1}$, and multiplying the inequalities (\ref{eq: inequalities Riccati sigma}) on the left and right by $Q_i$ yields
\[
Q_i A_i + A_i^T Q_i + Q_i W_i Q_i - \sum_{j=1}^M p_{ij} C^T_{ij} V^{-1}_{ij} C_{ij} \preceq 0,
\]
and hence we recover (\ref{eq: Riccati inequality}). In conclusion, the covariance matrices resulting from the switching policies approach within $O(\epsilon)$ as $\epsilon \to 0$ the covariance matrices which are obtained from the lower bound on the achievable cost. The theorem follows, by upper bounding the supremum limit of a sum by the sum of the supremum limits to get $Z^* \leq Z^\epsilon \leq Z^*+O(\epsilon)$.
\end{proof}

\begin{remark}
Since the bound computed in (\ref{eq: LMI relaxation}) is tight, and since it is easy to see that the performance bound of section \ref{section: 1D-dynamics} is at least as good as the bound (\ref{eq: LMI relaxation}) for the simplified problem of that section, we can conclude that the two bounds coincide and that section \ref{section: 1D-dynamics} gives an alternative way of computing the solution of  (\ref{eq: LMI relaxation}) in the case of identical sensors and one-dimensional systems. Using the closed form expression for the dual function (\ref{dual function - summary}), we only need to optimize over the unique Lagrange multiplier $\lambda$, independently of the number $N$ of systems, instead of solving the LMI (\ref{eq: LMI relaxation}), for which the number of variables grows with $N$.
\end{remark}

\subsubsection{Transient Behavior of the Switching Policies}

Before we conclude, we take a look at the transient behavior of the switching policies. We show that over a finite time interval, $\Sigma^\epsilon_i(t)$ remains close to $\tilde \Sigma_i(t)$, solution of the ``averaged" RDE (\ref{eq: averaged RDE}). Together with the previous result of lemma \ref{lem: convergence periodic} on the asymptotic behavior, we see then that $\Sigma^\epsilon_i(t)$ and $\tilde \Sigma_i(t)$ remain close for all $t$. For a matrix $A$, we denote by $\|A\|_{\infty}$ the maximal absolute value of the entries of $A$.

\begin{lemma}
For all $0 \leq T_0 < \infty$, there exist constants $\epsilon_0>0$ and $M_0 > 0$ such that for all $0<\epsilon\leq \epsilon_0$ and for all $t \in \left[ 0, T_0 \right]$, we have $\| \Sigma^\epsilon_i(t)-\tilde \Sigma_i(t) \|_\infty \leq M_0 \epsilon$.
\end{lemma}

\begin{proof}
As in the proof of lemma \ref{lem: approximation of average}, by Radon's lemma we know that $\Sigma_i^\epsilon$ is $C^\infty$ on each interval where $\sigma_i(t)$ is constant.
We have then, over the interval $t \in [l \epsilon, (l+\phi_1) \epsilon]$, for $l \in \mathbb{N}$:
\begin{equation}	\label{eq: DL}
\Sigma^\epsilon_i((l+\phi_1) \epsilon) = \Sigma^\epsilon_i({l \epsilon}) + \phi_1 \epsilon [A_i \Sigma^\epsilon_i({l \epsilon}) + \Sigma^\epsilon_i({l \epsilon}) A_i^T + W_i - \Sigma^\epsilon_i({l \epsilon}) C_{i j_1}^T V_{i j_1}^{-1} C_{i j_1} \Sigma^\epsilon_i({l \epsilon})] + O(\epsilon^2),
\end{equation}
where as before we denote $j_k:=\sigma_i(t)$ for $t \in [(l+\phi_1+\ldots+\phi_{k-1})\epsilon,(l+\phi_1+\ldots+\phi_{k})\epsilon]$.
Now over the period $t \in [(l+\phi_1) \epsilon, (l+\phi_1+\phi_{2}) \epsilon]$, we have:
\begin{align*}
\Sigma^\epsilon_i((l+\phi_1+\phi_{2}) \epsilon) &=\Sigma^\epsilon_i((l+\phi_1) \epsilon) + \phi_{2} \epsilon [A_i \Sigma^\epsilon_i((l+\phi_1) \epsilon)
+\Sigma^\epsilon_i((l+\phi_1) \epsilon)  A_i^T + W_i \\
&- \Sigma^\epsilon_i((l+\phi_1) \epsilon) C_{i j_{2}}^T V_{i j_{2}}^{-1} C_{i j_{2}} \Sigma^\epsilon_i((l+\phi_1) \epsilon) ] + O(\epsilon^2).
\end{align*}
Using (\ref{eq: DL}), we deduce that
\begin{align*}
\Sigma^\epsilon_i((l+\phi_1+\phi_{2}) \epsilon) = & \Sigma^\epsilon_i({l \epsilon}) + \epsilon [\phi_1+\phi_2] \{ A_i \Sigma^\epsilon_i({l \epsilon}) + \Sigma^\epsilon_i({l \epsilon}) A_i^T + W_i\} \\
& -\epsilon \Sigma^\epsilon_i({l \epsilon}) (\phi_1 C_{i j_1}^T V_{i j_1}^{-1} C_{i j_1} + \phi_2 C_{i j_{2}}^T V_{i j_{2}}^{-1} C_{i j_{2}})  \Sigma^\epsilon_i({l \epsilon}) + O(\epsilon^2).
\end{align*}
By immediate induction, and since $\phi_1+\cdots+\phi_K=1$, we then have
\[
\Sigma^\epsilon_i((l+1) \epsilon) =  \Sigma^\epsilon_i({l \epsilon}) + \epsilon \left \{A_i \Sigma^\epsilon_i({l \epsilon}) + \Sigma^\epsilon_i({l \epsilon}) A_i^T + W_i - \Sigma^\epsilon_i({l \epsilon}) \left( \sum_{k=1}^K \phi_k C_{i j_k}^T V_{i j_k}^{-1} C_{i j_k} \right)   \Sigma^\epsilon_i({l \epsilon}) \right \} + O(\epsilon^2).
\]
Hence by (\ref{eq: rewriting as frequencies}), $\Sigma^\epsilon_i$ verifies the relation
\begin{equation}	\label{eq: DL Sigma_epsilon}
\Sigma^\epsilon_i((l+1) \epsilon) =  \Sigma^\epsilon_i({l \epsilon}) + \epsilon \left \{A_i \Sigma^\epsilon_i({l \epsilon}) + \Sigma^\epsilon_i({l \epsilon}) A_i^T + W_i - \Sigma^\epsilon_i({l \epsilon}) \left( \sum_{j=1}^M p_{ij} C_{i j}^T V_{i j}^{-1} C_{i j} \right)   \Sigma^\epsilon_i({l \epsilon}) \right \} + O(\epsilon^2).
\end{equation}
But notice now that the approximation (\ref{eq: DL Sigma_epsilon}) is also true by definition for $\tilde \Sigma_i(t)$ over the interval $t \in [l \epsilon, (l+1)\epsilon]$. Next, consider the following identity for $Q, X$ and $\tilde X$ symmetric matrices:
\begin{align*}
&A \tilde X+ \tilde X A^T - \tilde X Q \tilde X - (A X+ X A^T - X Q X) \\
&=(A- \tilde X Q)(\tilde X-X) + (\tilde X-X)(A- \tilde X Q)^T + (\tilde X-X) Q (\tilde X-X).
\end{align*}
Letting $Q=\sum_{j=1}^M p_{ij} C_{i j}^T V_{i j}^{-1} C_{i j}$, $\Delta^\epsilon_i(l)=\tilde \Sigma_i(l \epsilon)-\Sigma^\epsilon_i(l \epsilon)$,
we obtain from this identity
\[
\Delta^\epsilon_i(l+1)=\Delta^\epsilon_i(l)+\epsilon \{ (A-\tilde \Sigma_i(l \epsilon) Q) \Delta^\epsilon_i(l) + \Delta^\epsilon_i(l) (A-\tilde \Sigma_i(l \epsilon) Q)^T + \Delta^\epsilon_i(l) Q \Delta^\epsilon_i(l) \} + O(\epsilon^2).
\]
Note that $\Delta^\epsilon_i(0)=0$ and $\tilde \Sigma_i(t)$ is bounded, so by immediate induction we have
\begin{flalign*}
&& & \Delta^\epsilon_i(l)=\sum_{k=1}^l R_k(\epsilon), \quad \text{ where } R_k(\epsilon)=O(\epsilon^2) \text{ for all } k. && \\
\text{Fix $T_0 \geq 0$. 
We have then}
&& & \tilde \Sigma_i \left( \left \lceil \frac{T_0}{\epsilon} \right \rceil \epsilon \right)-\Sigma^\epsilon_i \left( \left \lceil \frac{T_0}{\epsilon} \right \rceil \epsilon \right) = \Delta^\epsilon_i \left( \left \lceil \frac{T_0}{\epsilon} \right \rceil \right)=O(\epsilon). &&
\end{flalign*}
This means that there exist constants $\epsilon_0$, $M_0>0$ such that
\[
\left \| \tilde \Sigma_i \left( \left \lceil \frac{T_0}{\epsilon} \right \rceil \epsilon \right)-\Sigma^\epsilon_i \left( \left \lceil \frac{T_0}{\epsilon} \right \rceil \epsilon \right) \right \|_\infty \leq M_0 \epsilon, \quad \text{for all } 0 < \epsilon < \epsilon_0.
\]
It is easy to see from the argument above that a similar approximation is in fact valid for all $t$ up to time $\left \lceil \frac{T_0}{\epsilon} \right \rceil \epsilon$.
\end{proof}

\subsubsection{Numerical Simulation}

\begin{figure}
	\centering
	\includegraphics[width=0.7\textwidth]{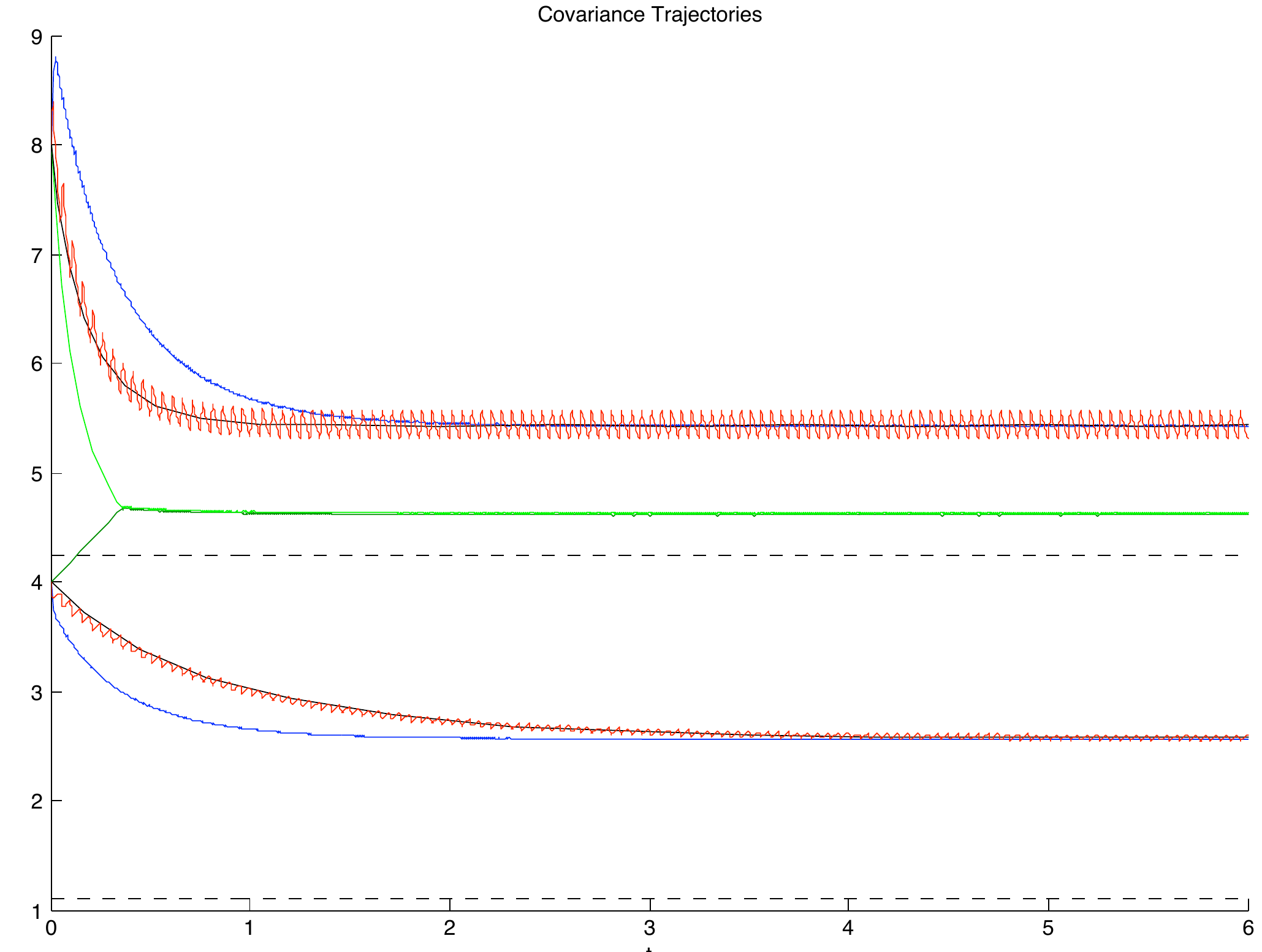}
	\caption{Comparison of the variance trajectories under Whittle's index policy (solid blue curves), the periodic switching policy (oscillating red curves), and the greedy policy (converging green curves). The solid black curves correspond to $\tilde \Sigma_i(t)$, solution of the RDE (\ref{eq: averaged RDE}). Here a single sensor switches between two scalar systems. The period $\epsilon$ was chosen to be $0.05$. The system parameters are $A_1=0.1, A_2=2, C_i=Q_i=R_i=1, \kappa_i=0$ for $i=1,2$. The dashed lines are the steady-state values that could be achieved with two identical sensors, each measuring one system. The performance of the Whittle policy is 7.98, which is optimal (i.e., matches the bound). The performance of the greedy policy is $9.2$. Note that the greedy policy makes the variances converge, while Whittle's policy makes the Whittle indices (not shown) converge. The switching policy spends $23\%$ of its time measuring system $1$ and $77\%$ of its time measuring system $2$.}
	\label{fig: simu}
\end{figure}

Figure \ref{fig: simu} compares the covariance trajectories for Whittle's index policy, the periodic switching policy and the greedy policy (measuring the system with highest mean square error on the estimate) for a simple problem with one sensor switching between two scalar systems. Significant improvements over the greedy policy can be obtained in general by using the periodic switching policies or the Whittle policy. An important computational advantage of the Whittle policy for large-scale problems with a limited number of identical sensors is that using the closed form solution of the indices provided in section \ref{section: summary Whittle}, it requires only ordering $N$ numbers (which is the same computational cost as for the greedy policy), whereas designing the open-loop switching policy requires computing a solution of the program (\ref{eq: LMI relaxation}).

\section{Conclusion}

In this paper, we have considered an attention-control problem in continuous time, which consists in scheduling sensor/target assignments and running the corresponding Kalman filters. We proved that the bound obtained from a RBP type relaxation is tight, assuming we allow the sensors to switch arbitrarily fast between the targets. An open question is to characterize the performance of the restless bandit index policy derived in the scalar case. It was found experimentally that the performance of this policy seems to match the bound, but we do not have a proof of this fact. An advantage of the Whittle index policy over the switching policies is that it is in feedback form. Obtaining optimal feedback policies for the multidimensional case would also be of interest. For practical applications, the main limitation of our model concerns the absence of switching costs and delays. Still, the optimal solution obtained in the absence of such costs should provide insight into the derivation of heuristics for more complex models. Additionally there are numerous sensor scheduling applications, such as for telemetry-data aerospace systems or radar waveform selection systems, where the switching costs are not too important.


\bibliography{biblio/biblio_thesis}

\end{document}